\documentclass[11pt]{amsart}

\usepackage{latexsym, graphics, psfrag, amscd, amssymb}
\newtheorem{thm}{Theorem}
\newtheorem{tm}[thm]{Theorem}
\newtheorem{fact}[thm]{Fact}

\newtheorem{rem}[thm]{Remark}

\newcommand{\ra}{\rightarrow}
\newcommand{\Si}{\Sigma}

\newcommand{\Z}{\mathbb{Z}}
\newcommand{\R}{\mathbb{R}}
\newcommand{\hol}{\mathrm{Hom}(\pi_1(\Si),G)}

\begin{document}

\parskip=0.35\baselineskip
\baselineskip=1.2\baselineskip

\title{Connected Components of The Space of Surface Group
Representations }

\author{Nan-Kuo Ho}
\address{Department of Mathematics\\ University of Toronto}
\email{nankuo@math.toronto.edu}

\author{Chiu-Chu Melissa Liu}
\address{Department of Mathematics \\Harvard University}
\email{ccliu@math.harvard.edu}


\subjclass{53}
\date{\today}

\begin{abstract}

Let $G$ be a connected, compact, semisimple Lie group. It is known
that for a compact closed orientable surface $\Sigma$ of genus $l
>1$, the order of the group $H^2(\Sigma,\pi_1(G))$ is equal to the
number of connected components of the space
$Hom(\pi_1(\Sigma),G)/G$ which can also be identified with the
moduli space of gauge equivalence classes of flat $G$-bundles over
$\Sigma$.
  We show that the same statement for a closed compact nonorientable surface which is homeomorphic to the connected sum of $k$ copies of the real projective plane, where $k\neq 1,2,4$, can be easily derived from a result in A. Alekseev, A.Malkin and E. Meinrenken's recent work on Lie group valued moment maps.

\end{abstract}

\maketitle

\section{introduction}

Let $\Si$ be a closed compact orientable surface of genus $l>1$.
In \cite{G1}, W.M. Goldman conjectured that for any connected complex
semisimple Lie group $G$, there is a bijection
\[
\pi_0(\hol)\ra H^2(\Si;\pi_1(G))\cong \pi_1(G).
\]
In \cite{Li}, J. Li proved the above conjecture
by Goldman as well as the following:

\begin{tm}[{\cite[Theorem 0.5]{Li}}]\label{orientable}
Let $\Si$ be a closed orientable Riemann surface of genus $l>1$.
Let $G$ be a connected, compact, semisimple Lie group. Then
there is a bijection
\[
\pi_0(\hol/G)\ra H^2(\Sigma;\pi_1(G))\cong \pi_1(G),
\]
where $G$ acts on $\hol$ by conjugation.
\end{tm}

Geometrically, $\hol/G$ can be identified with the moduli space of
gauge equivalence classes of flat $G$-bundles over $\Si$, where
a flat $G$-bundle is a principal $G$-bundle together with a flat
connection. It is known that there is a one-to-one correspondence
between topological principal $G$-bundles over $\Si$ and elements in
$H^2(\Si;\pi_1(G))\cong \pi_1(G)$.
In \cite{Li}, J. Li combined the argument in \cite{Ra} and \cite{AB}
to prove that the moduli space of gauge equivalence classes of flat
connections on a fixed underlying topological principal $G$-bundle
over $\Si$ is nonempty and connected, where $\Si$ and $G$ are as in
Theorem \ref{orientable}.

The main result in this paper is the following analogue of
Theorem \ref{orientable}.

\begin{tm}\label{nonorientable}
Let $\Si$ be a closed compact nonorientable surface which is
homeomorphic to the connected sum of $k$ copies of the real projective plane,
where $k\neq 1,2,4$.
Let $G$ be a connected, compact, semisimple Lie group. Then there is
a bijection
\[
\pi_0(\hol/G)\ra H^2(\Si;\pi_1(G))\cong \pi_1(G)/2\pi_1(G),
\]
where $2\pi_1(G)$ denote the subgroup
$\{k^2\mid k\in \pi_1(G)\}$ of the finite abelian group $\pi_1(G)$.
\end{tm}

One can see the isomorphism $H^2(\Si;\pi_1(G))\cong
\pi_1(G)/2\pi_1(G)$ as follows. By the universal coefficient
theorem of cohomology, we have a short exact sequence
\[
0 \ra \mathrm{Ext}(H_1(\Si;\Z),\pi_1(G))\ra H^2(\Si;\pi_1(G))
 \ra \mathrm{Hom}(H_2(\Si;\Z),\pi_1(G)) \ra 0,
\]
where $H_1(\Si,\Z)\cong \Z^{k-1}\oplus \Z/2\Z$, and $H_2(\Si,\Z)=0$.
So we have
\[
H^2(\Si,\pi_1(G))\cong \mathrm{Ext}(\Z/2\Z,\pi_1(G))\cong
\pi_1(G)/2\pi_1(G).
\]

Let $\Si$ be a closed compact nonorientable surface.
It is known that there is a one-to-one correspondence between
topological principal $G$-bundles over $\Si$ and
elements in $H^2(\Si;\pi_1(G))\cong \pi_1(G)/2\pi_1(G)$.
The geometric interpretation of Theorem \ref{nonorientable}
is that the moduli space of gauge equivalence classes of flat
connections on a fixed underlying topological principal $G$-bundle
over $\Si$ is nonempty and connected, where $\Si$ and $G$ are as in
Theorem \ref{nonorientable}. More details of this geometric
interpretation will be given in Appendix \ref{obstruction}.

We will show that both Theorem \ref{orientable}
and Theorem \ref{nonorientable} can be easily derived from the following
Fact \ref{fiber}, without identifying $\hol/G$ with the moduli
space of gauge equivalence classes of flat $G$-bundles over $\Si$.

\begin{fact}\label{fiber}
Let $G$ be a compact, semisimple, connected and simply connected Lie group. Let
$l$ be a positive integer. Then the commutator map
$\mu^l_G:G^{2l}\ra G$ defined by
\begin{equation}\label{commutator}
\mu^l_G(a_1,b_1,\ldots,a_l,b_l)= a_1b_1a_1^{-1}b_1^{-1}\cdots a_l
b_l a_l ^{-1}b_l^{-1}
\end{equation}
is surjective, and $(\mu^l_G)^{-1}(g)$ is connected for
all $g\in G$.
\end{fact}

The surjectivity follows from Goto's commutator theorem
\cite[Theorem 6.55]{HM}. The commutator map $\mu^l_G$ in Fact \ref{fiber}
is a group valued moment map of the q-Hamiltonian $G$-space
$G^{2l}$ in the sense of A. Alekseev, A. Malkin, and E. Meinrenken
\cite{AMM}. By \cite[Theorem 7.2]{AMM}, all the fibers of the moment map
of a connected q-Hamiltonian $G$-space are connected if $G$ is a
compact, connected and simply connected Lie group.

This paper is organized as follows. In Section \ref{riemann}, we
explain how Fact \ref{fiber} implies Theorem \ref{orientable} as
well as the genus one case. In Section \ref{odd} and Section
\ref{even}, we derive Theorem \ref{nonorientable} from Fact
\ref{fiber} for odd $k$ and even $k$, respectively.
In Section \ref{onetwofour}, we discuss the cases $k=1,2,4$
which are not covered by our main result.
In Appendix \ref{obstruction}, we give algebraic and geometric
interpretations of the obstruction map used in our proof.

\textit{Acknowledgment:} It is a pleasure to thank J.-K. Yu, who
motivated and guided this work. We also wish to thank L. Jeffrey
and E. Meinrenken for helpful conversations and valuable suggestions.

\section{Compact orientable surfaces}\label{riemann}

 In this section, $\Si$ is a closed Riemann surface with genus $l> 0$.
 Note that $\mathrm{Hom}(\pi_1(S^2),G)/G$ consists of a single point,
 so Theorem \ref{orientable} is not true for the genus zero case.

 Let $G$ be a connected, compact, semisimple Lie group, and let
 $\rho:\tilde{G}\ra G$ be the universal covering map which
 is also a group homomorphism. Then $\tilde{G}$ is a compact,
 semisimple, connected and simply connected Lie group, and $\mathrm{Ker}\rho$
 is a subgroup of the center $Z(\tilde{G})$ of $\tilde{G}$ which
 is abelian.

 Let $\mu_G^l:G^{2l}\ra G$ be the commutator map defined as
 in (\ref{commutator}), and let $e$ be the identity element of $G$.
 Then  $\hol$ can be identified with
\[
  (\mu_G^l)^{-1}(e)=\{(a_1,b_1,\ldots,a_l,b_l)\in G^{2l}\mid
  a_1 b_1 a_1^{-1}b_1^{-1}\cdots a_l b_l a_l^{-1} b_l^{-1}=e\}.
\]

There is an obstruction map (see e.g.\cite{G1})
\[
 o_2: \hol\cong (\mu_G^l)^{-1}(e) \ra \mathrm{Ker}\rho
\]
defined by
\[
(a_1,b_1,\ldots,a_l,b_l)\in (\mu_G^l)^{-1}(e)\mapsto
\tilde{a}_1\tilde{b}_1\tilde{a}_1^{-1}\tilde{b}_1^{-1}
\cdots\tilde{a}_l\tilde{b}_l\tilde{a}_l^{-1}\tilde{b}_l^{-1}
\in \mathrm{Ker}\rho,
\]
 where $(\tilde{a}_1,\tilde{b}_1\ldots,\tilde{a}_l,\tilde{b}_l)$
 is a preimage of $(a_1,b_1,\ldots,a_l,b_l)$ under
 $\rho^{2l}:\tilde{G}^{2l}\ra G^{2l}$. It is easily checked that the definition
 does not depend on the choice of
 $(\tilde{a}_1,\tilde{b}_1,\ldots,\tilde{a}_l,\tilde{b}_l)$
 and $o_2$ descends to a continuous map
\[
\bar{o}_2:\hol/G\ra\mathrm{Ker}\rho.
\]

For each $k\in\mathrm{Ker}\rho$,
there is a surjective continuous map
\[
(\mu^l_{\tilde{G}})^{-1}(k)\ra\bar{o}_2^{-1}(k)
\]
which is the restriction of
$\tilde{G}^{2l}\stackrel{\rho^{2l}}{\ra} G^{2l}\ra G^{2l}/G$,
where $G$ acts on $G^{2l}$ by diagonal conjugation. By Fact
\ref{fiber}, $(\mu^l_{\tilde{G}})^{-1}(k)$ is nonempty and
connected for each $k\in\mathrm{Ker}\rho$, so $\bar{o}_2$ is
surjective, and $\bar{o}_2^{-1}(k)$ is connected for each
$k\in\mathrm{Ker}\rho\cong\pi_1(G)$. This gives Theorem
\ref{orientable} as well as the genus one case.

\section{connected sum of a odd number of copies of the real projective plane}\label{odd}

In this section, $\Si$ is a closed compact nonorientable surface
which is homeomorphic to the connected sum of $2l+1$ copies of
$RP^2$, or equivalently, the connected sum of a Riemann surface of
genus $l$ and $RP^2$, where $l>0$.

Let $G$ be a connected, compact, semisimple Lie group, and let
$\rho:\tilde{G}\ra G$ be the universal covering map.
The space $\hol$ can be identified with
\[
X=\{ (a_1,b_1,\ldots, a_l,b_l,c)\in G^{2l+1} \mid
    a_1b_1a_1^{-1}b_1^{-1}\cdots a_l b_l a_l ^{-1}b_l^{-1} c^2=e \},
\]
where $e$ is the identity element of $G$.

Let $K=\mathrm{Ker}\rho\cong\pi_1(G)$, and
let $2K$ be the subgroup $\{k^2\mid k\in K\}$
of the finite abelian group $K$. There
is an obstruction map
\[
 o_2: \hol\cong X \ra K/2K
\]
defined by
\[
(a_1,b_1,\ldots,a_l,b_l,c)\in X\mapsto
 \tilde{a}_1 \tilde{b}_1\tilde{a}_1^{-1}\tilde{b}_1^{-1} \cdots
 \tilde{a}_l \tilde{b}_l \tilde{a}_l ^{-1}\tilde{b}_l^{-1} \tilde{c}^2
\in K/2K,
\]
 where $(\tilde{a}_1,\tilde{b}_1\ldots,\tilde{a}_l,\tilde{b}_l,\tilde{c})$
 is a preimage of $(a_1,b_1,\ldots,a_l,b_l,c)$ under
 $\rho^{2l+1}:\tilde{G}^{2l+1}\ra G^{2l+1}$.
 It is easily checked that the definition does not depend on the choice of
 $(\tilde{a}_1,\tilde{b}_1,\ldots,\tilde{a}_l,\tilde{b}_l,\tilde{c})$
 and $o_2$ descends to a continuous map
\[
\bar{o}_2:\hol/G\ra K/2K.
\]

For each $k\in K$, define
\[
\tilde{X}_k=\{
(a_1,b_1,\ldots,a_l,b_l,c)\in\tilde{G}^{2l+1}\mid
a_1 b_1 a_1^{-1}b_1^{-1}\cdots
a_l b_l a_l^{-1}b_l^{-1}c^2=k\}.
\]
Let $\bar{k}$ denote the image of $k$ under
$K\ra K/2K$. Then there is a surjective continuous map
\[
\tilde{X}_k\ra\bar{o}_2^{-1}(\bar{k})
\]
which is the restriction of
$\tilde{G}^{2l+1}\stackrel{\rho^{2l+1}}{\ra} G^{2l+1}\ra
G^{2l+1}/G$, where $G$ acts on $G^{2l+1}$ by diagonal conjugation.
Let $\tilde{e}$ denote the identity element of $\tilde{G}$. For
each $k\in K$, $\tilde{X}_k$ contains
$(\mu^l_{\tilde{G}})^{-1}(k)\times\{\tilde{e}\}$ which is nonempty
by Fact \ref{fiber}. Therefore, $\bar{o}_2$ is surjective. It
remains to show that $\tilde{X}_k$ is connected for all $k\in K$,
which will imply $\bar{o}_2^{-1}(\bar{k})$ is
connected for all $\bar{k}\in K/2K$.

Let $Q': \tilde{G}^{2l+1}\ra \tilde{G}$ be the projection to the
last factor. The restriction of $Q'$ to $\tilde{X}_k$ gives a map
$Q:\tilde{X}_k\ra \tilde{G}$,
$(a_1,b_1,\ldots,a_l,b_l,c)
\mapsto c$. Fact \ref{fiber} implies that $Q$ is surjective and
$Q^{-1}(c)$ is connected for all $c\in \tilde{G}$.
Let $q\in T$ be a square root of $k\in Z(\tilde{G})\subset T$,
where we fix a maximal torus $T$ of $\tilde{G}$.
To prove that $\tilde{X}_k$ is connected, it
suffices to show that for any $c\in \tilde{G}$, there is a path
$\gamma:[0,1]\ra \tilde{X}_k$ such that
$\gamma(0)\in Q^{-1}(q)$ and
$\gamma(1)\in Q^{-1}(c)$.

For any $c\in \tilde{G}$, there exists $g\in \tilde{G}$
such that $g^{-1}c g\in T$.
Let $\tilde{\mathfrak{g}}$ and $\mathfrak{t}$ be the Lie algebras of
$\tilde{G}$ and $T$, respectively.
Let $\exp:\tilde{\mathfrak{g}}\ra\tilde{G}$ be the exponential map.
Then $g^{-1}cgq^{-1}=\exp\xi$ for some $\xi\in\mathfrak{t}$.
Let $W$ be the Weyl group of $\tilde{G}$, and let $w\in W$ be a
Coxeter element. The linear map $w:\mathfrak{t}\ra\mathfrak{t}$
has no eigenvalue equal to $1$ \cite[Section 3.16]{Hu},
so there exists $\xi'\in\mathfrak{t}$ such that $w\cdot \xi'-\xi'=\xi$.
Recall that $W=N(T)/T$, where $N(T)$ is the normalizer of
$T$ in $\tilde{G}$, so $w=aT\in N(T)/T$ for some $a\in \tilde{G}$.
We have
\[
a\exp(t\xi')a^{-1}\exp(-t\xi')=\exp(t\xi)
\]
for any $t\in\R$.

The group $\tilde{G}$ is connected, so there exists a path
$\tilde{g}:[0,1]\ra\tilde{G}$ such that
$\tilde{g}(0)=\tilde{e}$ and $\tilde{g}(1)=g$.
Define $\gamma:[0,1]\ra \tilde{G}^{2l+1}$ by
\[
\gamma(t)=(a(t),b(t),\tilde{e},\ldots,\tilde{e},c(t)),
\]
where
\[
a(t)=\tilde{g}(t)a \tilde{g}(t)^{-1},\
b(t)=\tilde{g}(t)\exp(-2t\xi')\tilde{g}(t)^{-1},\
c(t)=\tilde{g}(t)q\exp(t\xi)\tilde{g}(t)^{-1}.
\]
Then the image of $\gamma$ lies in $\tilde{X}_k$,
$\gamma(0)=(a,\tilde{e},\tilde{e},\ldots,\tilde{e},q)\in
Q^{-1}(q)$, and
\[
\gamma(1)=(gag^{-1},g\exp(-2\xi')g^{-1},
\tilde{e},\ldots,\tilde{e},c)\in Q^{-1}(c).
\]

\begin{rem}
Let $M$ be a manifold, and let $G$ be a compact Lie group.
By \cite[Definition 2.2 (B1)]{AMM}, a necessary condition
for a map $\mu: M\ra G$ to be a group valued moment map is
$\mu^*\chi$ being exact, where $\chi$ is the canonical closed
bi-invariant $3$-form on $G$.

Let $\tilde{G}$ be as above, and define
$\tilde{\mu}: \tilde{G}^{2l+1} \ra \tilde{G}$ by
\[
(a_1,b_1,\ldots,a_l,b_l,c) \mapsto
a_1b_1a_1^{-1}b_1^{-1}\cdots a_l b_l a_l ^{-1}b_l^{-1} c^2.
\]
Let $Q_{2l+1}:\tilde{G}^{2l+1}\ra\tilde{G}$ be the projection
to the last factor. We have $\tilde{\mu}^*\chi=2Q_{2l+1}^*\chi$,
which represents a nontrivial cohomology class. So
$\tilde{\mu}$ cannot be a group valued moment map,
and we cannot apply \cite[Theorem 7.2]{AMM} directly
to $\tilde{\mu}$ to conclude that $\tilde{X}_k=\tilde{\mu}^{-1}(k)$ is connected.
\end{rem}

\section{connected sum of an even number of copies of the real projective plane}\label{even}

In this section, $\Si$ is a compact nonorientable surface which is
homeomorphic to the connected sum of $2l+2$ copies of $RP^2$, or
equivalently, the connected sum of a Riemann surface of genus $l$
and a Klein bottle, where $l>1$.

Let $G$ be a connected, compact, semisimple Lie group, and let
$\rho:\tilde{G}\ra G$ be the universal covering map.
The space $\hol$ can be identified with
\[
X=\{ (a_1,b_1,\ldots, a_l,b_l,c_1,c_2)\in G^{2l+2} \mid
    a_1b_1a_1^{-1}b_1^{-1}\cdots a_l b_l a_l ^{-1}b_l^{-1} c_1^2 c_2^2=e \},
\]
where $e$ is the identity element of $G$.

Let $K=\mathrm{Ker}\rho\cong\pi_1(G)$.
There is an obstruction map
\[
 o_2: \hol\cong X \ra K/2K
\]
 defined by
\[
(a_1,b_1,\ldots,a_l,b_l,c_1,c_2)\in X\mapsto
 \tilde{a}_1 \tilde{b}_1\tilde{a}_1^{-1}\tilde{b}_1^{-1} \cdots
 \tilde{a}_l \tilde{b}_l \tilde{a}_l ^{-1}\tilde{b}_l^{-1}
 \tilde{c}_1^2\tilde{c}_2^2 \in K/2K,
\]
 where $(\tilde{a}_1,\tilde{b}_1\ldots,\tilde{a}_l,\tilde{b}_l,
 \tilde{c}_1,\tilde{c}_2)$
 is a preimage of $(a_1,b_1,\ldots,a_l,b_l,c_1,c_2)$ under
 $\rho^{2l+2}:\tilde{G}^{2l+2}\ra G^{2l+2}$.
 It is easily checked that the definition does not depend on the choice of
 $(\tilde{a}_1,\tilde{b}_1,\ldots,\tilde{a}_l,\tilde{b}_l,
 \tilde{c}_1,\tilde{c}_2)$
 and $o_2$ decends to a continuous map
\[
\bar{o}_2:\hol/G\ra K/2K.
\]

For each $k\in K$, define
\[
\tilde{X}_k=\{
(a_1,b_1,\ldots,a_l,b_l,c_1,c_2)\in\tilde{G}^{2l+2}\mid
a_1 b_1 a_1^{-1}b_1^{-1}\cdots a_l b_l a_l^{-1}b_l^{-1}
c_1^2 c_2^2=k\}.
\]
Let $\bar{k}$ be the image of $k$ under $K\ra K/2K$.
Then there is a surjective continuous map
\[
\tilde{X}_k\ra\bar{o}_2^{-1}(\bar{k})
\]
which is the restriction of
$\tilde{G}^{2l+2}\stackrel{\rho^{2l+2}}{\ra} G^{2l+2}\ra
G^{2l+2}/G$, where $G$ acts on $\tilde{G}^{2l+2}$ by diagonal
conjugation. For each $k\in K$, $\tilde{X}_k$ contains
$(\mu^l_{\tilde{G}})^{-1}(k)\times\{e\}\times\{e\}$ which is
nonempty by Fact \ref{fiber}. Therefore, $\bar{o}_2$ is
surjective. It remains to show that $\tilde{X}_k$ is connected for
all $k\in K$, which will imply that
$\bar{o}_2^{-1}(\bar{k})$ is connected for all
$\bar{k}\in K/2K$.

Let $Q': \tilde{G}^{2l+2}\ra \tilde{G}^2$ be the projection to the
last two factors. The restriction of $Q'$ to $\tilde{X}_k$ gives a map
$Q:\tilde{X}_k\ra \tilde{G}^2$,
$(a_1,b_1,\ldots,a_l,b_l,c_1,c_2)\mapsto (c_1,c_2)$.
Fact \ref{fiber} implies that $Q$ is surjective and
$Q^{-1}(c_1,c_2)$ is connected for all
$(c_1,c_2)\in \tilde{G}^2$.
Let $q\in T$ be a square root of $k\in Z(\tilde{G})$, where
we fix a maximal torus $T$ of $\tilde{G}$.
To prove that $\tilde{X}_k$ is connected, it
suffices to show that for any $(c_1,c_2)\in \tilde{G}^2$,
there is a path $\gamma:[0,1]\ra \tilde{X}_k$ such that
$\gamma(0)\in Q^{-1}(\tilde{e},q)$ and
$\gamma(1)\in Q^{-1}(c_1,c_2)$,
where $\tilde{e}$ is the identity element of $\tilde{G}$.

For any $(c_1,c_2)\in \tilde{G}^2$, there exists
$g_1,g_2\in \tilde{G}$ such that $g_j^{-1}c_j g_j\in T$ for $j=1,2$.
Let $\tilde{\mathfrak{g}}$ and $\mathfrak{t}$ be the Lie algebras of
$\tilde{G}$ and $T$, repsectively.
Let $\exp:\tilde{\mathfrak{g}}\ra\tilde{G}$ be the exponential map.
Then $g_1^{-1}c_1 g_1=\exp\xi_1$ for some $\xi_1\in\mathfrak{t}$, and
$g_2^{-1}c_2 g_2q^{-1}=\exp\xi_2$ for some $\xi_2 \in\mathfrak{t}$.
Let $W$ be the Weyl group of $\tilde{G}$, and let $w\in  W$ be a Coxeter
element. The linear map $w:\mathfrak{t}\ra\mathfrak{t}$ has no
eigenvalue equal to $1$, so there exists
$\xi_j'\in\mathfrak{t}$ such that $w\cdot \xi_j'-\xi_j'=\xi_j$
for $j=1,2$, where $w=a T\in N(T)/T=W$ for some $a\in \tilde{G}$.
We have
\[
a\exp(t\xi_j')a^{-1}\exp(-t\xi_j')=\exp(t\xi_j)
\]
for any $t\in\R$, where $j=1,2$.

The group $\tilde{G}$ is connected, so there exist paths
$\tilde{g}_j:[0,1]\ra\tilde{G}$ such that
$\tilde{g}_j(0)=\tilde{e}$ and $\tilde{g}_j(1)=g_j$ for $j=1,2$.
Define $\gamma:[0,1]\ra \tilde{G}^{2l+2}$ by
\[
\gamma(t)=(a_2(t),b_2(t),a_1(t),b_1(t),\tilde{e},\ldots,\tilde{e},
c_1(t),c_2(t))
\]
where
\[
a_j(t)=\tilde{g}_j(t)a \tilde{g}_j(t)^{-1},\
b_j(t)=\tilde{g}_j(t)\exp(-2t\xi'_j )\tilde{g}_j(t)^{-1}
\]
for $j=1,2$, and
\[
c_1(t)=\tilde{g}_1(t) \exp(t\xi_1)\tilde{g}_1(t)^{-1},\
c_2(t)=\tilde{g}_2(t)q\exp(t\xi_2)\tilde{g}_2(t)^{-1}.
\]
Then the image of $\gamma$ lies in $\tilde{X}_k$,
$\gamma(0)=(a,\tilde{e},a,\tilde{e},\tilde{e},\ldots,
\tilde{e},q)\in Q^{-1}(\tilde{e},q)$,
and
\[
\gamma(1)=(g_2 a g_2^{-1}, g_2\exp(-2\xi_2')g_2^{-1},
g_1 a g_1^{-1},g_1\exp(-2\xi_1')g_1^{-1},
\tilde{e},\ldots,\tilde{e},c_1,c_2)\in Q^{-1}(c_1,c_2).
\]

\begin{rem}

Let $\tilde{G}$ be as above, and define
$\tilde{\mu}: \tilde{G}^{2l+2} \ra \tilde{G}$ by
\[
(a_1,b_1,\ldots,a_l,b_l,c_1,c_2) \mapsto
a_1b_1a_1^{-1}b_1^{-1}\cdots a_l b_l a_l ^{-1}b_l^{-1} c_1^2 c_2^2.
\]
Let $Q_{2l+1}$ and $Q_{2l+2}$ denote the projections
from $G^{2l+2}$ to the $(2l+1)$-th and $(2l+2)$-th
factors, respectively. We have
$\tilde{\mu}^*\chi=2Q_{2l+1}^*\chi+ 2Q_{2l+2}^*\chi$,
which represents a nontrivial cohomology class. So
$\tilde{\mu}$ cannot be a group valued moment map,
and we cannot apply \cite[Theorem 7.2]{AMM} directly to
$\tilde{\mu}$ to conclude that $\tilde{X}_k=\tilde{\mu}^{-1}(k)$ is connected.
\end{rem}

\section{The cases $k=1,2,4$}\label{onetwofour}

Let $G$ be a compact, semisimple, connected and simply connected
Lie group. By Theorem \ref{nonorientable}, $\hol/G$ is connected
if $\Si$ is the connected sum of $k$ copies of $RP^2$, where
$k\neq 1,2,4$. For $k=1$, we have
\[
\mathrm{Hom}(\pi_1(RP^2),G)/G=
\{ g\in G \mid g^2 =e\}/G,
\]
where $G$ acts by conjugation. In particular, for $n>1$,
\[
\mathrm{Hom}(\pi_1(RP^2),SU(n))/SU(n)\cong
\left\{\left.\left(\begin{array}{cc}-I_{2j}&0\\0& I_{n-2j}\end{array}
\right)\right|
j=0,1,\ldots,\left[\frac{n}{2}\right]\right\},
\]
where $I_m$ denotes the $m\times m$ identity matrix. For $n>2$,
\[
\mathrm{Hom}(\pi_1(RP^2),Spin(n))/Spin(n)\cong
\left\{ (-1)^j e_1 e_2\cdots e_{2j}\mid j=0,1,\ldots,\left[\frac{n}{2}\right]
\right\},
\]
where we view $Spin(n)$ as a subset of the real Clifford algebra $C_n$,
and $\{e_1,\ldots,e_n\}$ is an orthonormal basis for $\R^n$. For $n>0$,
\[
\mathrm{Hom}(\pi_1(RP^2),Sp(n))/Sp(n)\cong\left\{
\left.\left(\begin{array}{cc}
\begin{array}{cc}-I_k&0\\0&I_{n-k}\end{array} &0\\
0&\begin{array}{cc}-I_k&0\\0&I_{n-k}\end{array}
\end{array}\right)\right|k=0,1,\ldots,n
\right\}.
\]
So Theorem \ref{nonorientable} is not true for $k=1$.

It is not clear to us if Theorem \ref{nonorientable} is
true for $k=2$. We expect Theorem \ref{nonorientable}
to hold for $k=4$, but we do not know how to prove it using
the approach in this paper.

\appendix

\section{Obstruction class} \label{obstruction}
Let $\Si$ be a connected closed compact surface.
Let $G$ be a connected, compact, semisimple Lie group.
Let $\rho:\tilde{G}\ra G$ be the universal covering map which
is also a group homomorphism. Then $\mathrm{Ker}\rho$
is contained in the center $Z(\tilde{G})$ of $\tilde{G}$.
Let $K=\mathrm{Ker\rho}\cong\pi_1(G)$ which is a finite abelian
group.

In this appendix, we will give several interpretations
of the obstruction map
\[
o_2:\hol\ra H^2(\Si;K).
\]
Note that the definition in Section \ref{odd} can be extended
to the case $l=0$, and the definition in Secton \ref{even} can be
extended to the cases $l=0, 1$.
If $\Si$ is homeomorphic to a 2-sphere, then $\hol$ consists of
a single point $e$. We extend the definition in Section
\ref{riemann} to the genus zero case by defining $o_2(e)$ to be the
identity element of $H^2(\Si;K)$. So the obstruction map $o_2$ is defined
for any closed compact surfaces $\Si$.

\subsection{Non-abelian cohomology of groups}

For an account of non-abelian cohomology of groups, see
\cite[Chapter VII, Appendix]{Se}.

We have a short exact sequence of groups
\begin{equation}\label{exact}
1\ra K\ra \tilde{G} \ra G\ra 1.
\end{equation}

Let $\pi_1(\Si)$ act trivially on (\ref{exact}) on
the left so that (\ref{exact}) can be viewed as a short exact sequence of
(non-abelian) $\pi_1(\Si)$-modules.
The short exact sequence (\ref{exact}) gives rise to the following long exact
sequence of pointed sets
\cite[Chapter VII, Appendix, Proposition 2]{Se}:
\begin{equation}\label{group}
\begin{array}{ccccccccc}
1 &\ra& H^0(\pi_1(\Si),K) &\ra& H^0(\pi_1(\Si),\tilde{G}) &\ra&
    H^0(\pi_1(\Si),G)&&\\
  &\stackrel{\delta}{\ra}&
   H^1(\pi_1(\Si),K) &\ra& H^1(\pi_1(\Si),\tilde{G}) &\ra&
    H^1(\pi_1(\Si),G) &\stackrel{\Delta}{\ra}& H^2(\pi_1(\Si),K)
\end{array}
\end{equation}
which can be rewritten as
\begin{equation}\label{hom}
\begin{array}{ccccccccc}
1&\ra& K&\ra&\tilde{G}&\ra& G&&\\
&\stackrel{\delta}{\ra}& \mathrm{Hom}(\pi_1(\Si),K)&\ra&
 \mathrm{Hom}(\pi_1(\Si),\tilde{G})/\tilde{G}&\ra& \hol/G
&\stackrel{\Delta}{\ra}& H^2(\pi_1(\Si),K)
\end{array}
\end{equation}
by checking the definitions.

Let $p\in\hol$ represent a point $\bar{p}\in H^1(\pi_1(\Si),G)\cong\hol/G$.
From the exact sequence of (\ref{hom}),
$p$ can be lifted to a homomophism $\tilde{p}:\pi_1(\Si)\ra \tilde{G}$
if and only if $\Delta(\bar{p})$ is the identity element of
$H^2(\pi_1(\Si),K)$. From the explicit expression of the obstruction map
$o_2:\hol\ra H^2(\Si;K)$ we see that
$o_2(p)$ is the identity element of the group $H^2(\Si;K)$ if and only
if $p$ can be lifted to a homomorphism $\tilde{p}:\pi_1(\Si)\ra\tilde{G}$.
Actually,
\[
\Delta: H^1(\pi_1(\Si), G)\ra H^2(\pi_1(\Si),K)
\]
coincides with
\[
\bar{o}_2:\hol/G\ra H^2(\Si;K)
\]
under the identifications
\[
H^1(\pi_1(\Si), G)\cong \hol/G,\ H^2(\pi_1(\Si),K)\cong H^2(\Si;K),
\]
where $H^2(\pi_1(\Si),K)\cong H^2(\Si;K)$ because $\Si$ is the
Eilenberg-MacLane space $K(\pi_1(\Si),1)$.

\subsection{Non-abelian \v{C}ech cohomology} For an account of
of non-abelian \v{C}ech cohomology, see \cite[Appendix A]{LM}.

\subsubsection{Flat bundles}
Given a group $A$, let $A$ also denote the sheaf on $\Si$
for which $A(U)$ is the group of locally constant functions from $U$ to $A$,
where $U$ is any open subset of $\Si$. The short exact sequence of sheaves
\[
1\ra K\ra \tilde{G}\ra G\ra 1
\]
gives rise to the following long exact sequences
of non-abelian \v{C}ech cohomology \cite[Remark A.2]{LM}
\begin{equation}\label{cech}
\begin{array}{ccccccccc}
1 &\ra& H^0(\Si;K) &\ra& H^0(\Si;\tilde{G}) &\ra&
    H^0(\Si;G)&&\\
  &\stackrel{\delta}{\ra}&
   H^1(\Si;K) &\ra& H^1(\Si;\tilde{G}) &\ra&
    H^1(\Si;G) &\stackrel{\Delta}{\ra}& H^2(\Si;K)
\end{array}
\end{equation}
which can be rewritten as
\begin{equation}\label{hom2}
\begin{array}{ccccccccc}
1&\ra& K&\ra&\tilde{G}&\ra& G&&\\
&\stackrel{\delta}{\ra}& H^1(\Si;K)&\ra&
 H^1(\Si;\tilde{G})&\ra& H^1(\Si;G)
&\stackrel{\Delta}{\ra}& H^2(\Si;K)
\end{array}
\end{equation}
by checking the definitions.
The pointed sets $H^1(\Si;\tilde{G})$ and
$H^1(\Si;G)$ can be identified with the moduli
spaces of gauge equivalance classes of flat $\tilde{G}$-bundles
and $G$-bundles, respectively. The map $\delta$ in (\ref{hom2})
is the trivial group homomorphism, while the map $\Delta$ in
(\ref{hom2}) sends a gauge equivalence class of flat $G$-bundles
to the obstruction to lifting it to a gauge equivalence
class of flat $\tilde{G}$-bundles. Actually,
\[
\Delta: H^1(\Si;G)\ra H^2(\Si;K)
\]
coincides with
\[
\bar{o}_2: \hol/G\ra H^2(\Si;K)
\]
under the identification $H^1(\Si;G)\cong \hol/G$.

\subsubsection{Topological bundles} \label{topbundle}
Given a group $A$, let $\underline{A}$ denote the sheaf on $\Si$
for which $\underline{A}(U)$ is the group of continuous functions
from $U$ to $A$, where $U$ is any open subset of $\Si$.
Note that $\underline{K}=K$ since $K$ is discrete.
The short exact sequence of sheaves
\[
1\ra K\ra \underline{\tilde{G}}\ra \underline{G}\ra 1
\]
gives rise to the following long exact sequence of \v{C}ech
cohomology
\begin{equation}\label{topology}
\begin{array}{ccccccccc}
1 &\ra& H^0(\Si;K) &\ra& H^0(\Si;\underline{\tilde{G}}) &\ra&
    H^0(\Si,\underline{G})&&\\
  &\stackrel{\delta'}{\ra}&
   H^1(\Si;K) &\ra& H^1(\Si;\underline{\tilde{G}}) &\ra&
    H^1(\Si;\underline{G}) &\stackrel{\Delta'}{\ra}& H^2(\Si;K).
\end{array}
\end{equation}
The pointed sets $H^1(\Si;\underline{\tilde{G}})$ and
$H^1(\Si;\underline{G})$ can be identified with
the set of equivalence classes of topological principal
$\tilde{G}$-bundles and $G$-bundles, respectively.
Note that $H^1(\Si;K)\cong \mathrm{Hom}(\pi_1(\Si),K)$.
The exact sequence (\ref{topology}) can be rewritten as
\begin{equation}
\begin{array}{ccccccccc}
1&\ra&K&\ra&\mathrm{Map}(\Si,\tilde{G})&\ra&
\mathrm{Map}(\Si,G)&&\\
&\stackrel{\delta'}{\ra}& \mathrm{Hom}(\pi_1(\Si),K)&\ra&
 \mathrm{Prin}_{\tilde{G}}(\Si)&\ra& \mathrm{Prin}_G(\Si)
&\stackrel{\Delta'}{\ra}& H^2(\Si,K).
\end{array}
\end{equation}
The map $\delta'$ sends $f:\Si\ra G$ to
$f_*:\pi_1(\Si)\ra \pi_1(G)\cong K$,
and $\Delta'$ sends a principal $G$-bundle to the obstruction to lifting
it to a principal $\tilde{G}$-bundle. For example, if $G=SO(n)$ ($n>2$), then
$K=\Z/2\Z$, and $\Delta'(P)$ is the second Stiefel-Whitney
class $w_2(P)\in H^2(\Si,\Z/2\Z)$.

We have a commutative diagram of sheaves
\[
\begin{CD}
1 @>>> K @>>> \tilde{G} @>>> G @>>> 1\\
&&  @VVV @VVV @VVV\\
1 @>>> K @>>> \underline{\tilde{G}} @>>> \underline{G} @>>> 1
\end{CD}
\]
which induces a morphism from (\ref{cech}) to (\ref{topology}). In
particular, we have the following commutative diagram
\[
\begin{CD}
H^1(\Si;G)=&\hol/G@>{\Delta}>> H^2(\Si;K)\\
@V{i}VV & @V{id}VV\\
H^1(\Si;\underline{G})=&\mathrm{Prin}_G(\Si) @>{\Delta'}>> H^2(\Si;K)
\end{CD}
\]
where $i$ sends a flat $G$-bundle to the underlying topological principal
$G$-bundle. Therefore, if $P$ is the underlying topological principal
$G$-bundle of a flat $G$-bundle associated to $p\in \hol$, then
$\Delta'(P)=o_2(p)\in H^2(\Si;K)$.

\subsection{Obstruction theory} For an account of obstruction theory
of fiber bundles, see \cite{St}.

The first obstruction to the triviality of a principal $G$-bundle $P$
over $\Si$, or equivalently, the first obstruction of existence of a
cross-section, lies in $H^1(\Si;\pi_0(G))$, which vanishes because $G$
is connected. So $P$ is trival over the 1-skeleton.
The second obstruction lies in $H^2(\Si;\pi_1(G))=H^2(\Si;K)$
and coincides with $\Delta'(P)$, the obstruction to lifting $P$ to a
principal $\tilde{G}$-bundle. For a surface $\Si$, $\Delta'(P)$ is the
only obstruction to the triviality of $P$. So the topological principal
$G$-bundles over $\Si$ are classified by $\Delta'(P)\in H^2(\Si;K)$.
The set $\mathrm{Prin}_{\tilde{G}}(\Si)$ consists of a single point
corresponding to the trivial $\tilde{G}$-bundle, and
$\Delta'$ in (\ref{topology}) is a bijection.

Note that if $G$ is not simply connected, then there exists nontrivial
topological $G$-bundles over $S^2$, but the underlying topological principal
$G$-bundle of a flat $G$-bundle over $S^2$ must be trivial.
So $\Delta$ sends the only point in
$H^1(S^2;G)=H^1(\pi_1(S^2),G)$) to the identity
element of $H^2(S^2;K)=H^2(\pi_1(S^2),K)$.

\subsection{Elementary approach}
There is an elementary approach to the classification of topological
principal $G$-bundles over $\Si$ which we learned from
E. Meinrenken. Let $D$ be a disc around a point $x\in \Si$. Any principal
$G$-bundle $P$ is trival over $D$ and over $\Si\setminus\{x\}$,
and the topological type of $P$ is determined by the homotopy class of
the transition function  $\psi: D\setminus\{x\}\ra G$, or equivalently,
an element in $K=\pi_1(G)$.
This gives a surjective map $K\ra \mathrm{Prin}_G(\Si)$
which is injective if $\Si$ is orientable
and induces a bijection $K/2K\ra \mathrm{Prin}_G(\Si)$
if $\Si$ is nonorientable. Therefore, we have a bijection
\[
 \phi: H^2(\Si;K)\ra \mathrm{Prin}_G(\Si).
\]
The argument in \cite[Section 6]{Li} shows that if $P$ is the underlying
topological principal $G$-bundle of a flat $G$-bundle associated to
$p\in \hol$, then $\phi(o_2(p))=P$.

Based on above discussions, Theorem \ref{orientable} and
Theorem \ref{nonorientable} can be reformulated as follows.
\begin{tm}
Let $\Si$ be a closed compact orientable surface of genus $l>0$,
or a closed compact nonorientable surface which is
homeomorphic to $k$ copies of the real projective plane, where
$k\neq 1,2,4$. Let $G$ be a connected, compact, semisimple
Lie group, and let $P$ be a topological principal $G$-bundle
over $\Si$. Then the moduli space of gauge equivalence classes
of flat connections on $P$ is nonempty and connected.
\end{tm}


\begin{thebibliography}{AMW}

\bibitem[AB]{AB} M.F. Atiyah and R. Bott,
{\em The Yang-Mills equations over Riemann surfaces},
Philos. Trans. Roy. Soc. London Ser. A \textbf{308} (1983),
no. 1505, 523--615.

\bibitem[AMM]{AMM} A. Alekseev, A. Malkin, and E. Meinrenken,
 {\em Lie group valued moment maps},
 J. Differential Geom. \textbf{48} (1998), no. 3, 445--495.

\bibitem[AMW]{AMW} A. Alekseev, E. Meinrenken, and C. Woodward,
 {\em Duistermaat-Heckman measures and moduli spaces of flat bundles over
 surfaces},
 Geom. Funct. Anal. \textbf{12} (2002), no. 1, 1--31.

\bibitem[BD]{bd} Theodor Br\"{o}cker and Tammo tom Dieck,
 {\em Representations of compact Lie groups},
 Graduate Texts in Mathematics, 98. Springer-Verlag, New York, 1985.

\bibitem[Go]{G1} W.M. Goldman,
 {\em Topological components of spaces of representations},
 Invent. Math. \textbf{93} (1988), no. 3, 557--607.

\bibitem[GM]{GM} W.M. Goldman and J.J. Millson,
{\em The deformation theory of representations of fundamental groups
of compact K\"{a}hler manifolds},
Bull. Amer. Math. Soc. (N.S.) \textbf{18} (1988),
no. 2, 153--158.

\bibitem[HM]{HM} K.H. Hofmann and S.A. Morris,
{\em The structure of compact groups.
A primer for the student---a handbook for the expert},
de Gruyter Studies in Mathematics, 25.
Walter de Gruyter \& Co., Berlin, 1998.

\bibitem[Hu]{Hu} J.E. Humphreys,
{\em Reflection groups and Coxeter groups},
Cambridge Studies in Advanced Mathematics, 29.
Cambridge University Press, Cambridge, 1990.

\bibitem[Li]{Li} J. Li,
 {\em The space of surface group representations},
 Manuscripta Math. \textbf{78} (1993), no. 3, 223--243.

\bibitem[LM]{LM} H.B. Lawson and M. Michelsohn,
{\em Spin geometry},
Princeton Mathematical Series 38,
Princeton University Press, Princeton, NJ, 1989.

\bibitem[Mi]{Mi} J. Milnor,
{\em On the existence of a connection with curvature zero},
Comment. Math. Helv. \textbf{32} (1958), 215--223.

\bibitem[Ra]{Ra} A. Ramanathan,
{\em Moduli for principal bundles},
Algebraic geometry (Proc. Summer Meeting,
Univ. Copenhagen, Copenhagen, 1978), pp. 527--533, Lecture Notes in Math.,
732, Springer, Berlin, 1979.

\bibitem[Se]{Se} J.-P. Serre,
{\em Local fields},
translated from the French by Marvin Jay Greenberg,
Graduate Texts in Mathematics, 67. Springer-Verlag, New York-Berlin, 1979.

\bibitem[St]{St} N. Steenrod,
{\em The Topology of Fibre Bundles},
Princeton Mathematical Series, vol. 14. Princeton University Press, Princeton, N. J., 1951.

\end{thebibliography}
\end{document}